\documentclass[11pt,reqno]{amsart}
\usepackage{amsmath,amssymb,latexsym, url,graphicx,newtxtext,newtxmath,tikz}

\usepackage[margin=1.15in]{geometry}
\usepackage{booktabs}
\usepackage{xcolor}
\usepackage{enumitem}
\usepackage[colorlinks=true,linkcolor=blue!50!black,citecolor=blue!50!black,urlcolor=blue!50!black]{hyperref}

\newtheorem{theorem}{Theorem}[section]
\newtheorem{lemma}[theorem]{Lemma}

\theoremstyle{definition}

\theoremstyle{remark}

\newcommand{\dg}{\,d\gamma}

\begin{document}

\title{An improved lower bound for piercing the chessboard by lines}

\author{Ivan Mati\'{c}}
\address{Department of Mathematics, Baruch College, City University of New York,
One Bernard Baruch Way, New York, NY 10010, USA}
\email{ivan.matic@baruch.cuny.edu}

\author{Rado\v{s} Radoi\v{c}i\'{c}}
\address{Department of Mathematics, Baruch College, City University of New York,
One Bernard Baruch Way, New York, NY 10010, USA}
\email{rados.radoicic@baruch.cuny.edu}

\date{\today}

\begin{abstract}
Let $p_n$ be the least number of lines meeting the interior of every cell of an $n\times n$ board.
We prove that $\liminf_{n\to\infty}p_n/n\ge0.7520517>\frac34$. The previous bound was $0.7$.
The proof is computer-assisted. It identifies a non-negative polynomial on the unit square and bounds its integral along every line.
\end{abstract}

\maketitle

\section{Setup and results}\label{sec:results}

Let $B_n=[0,n]^2$. Its cells are the open unit squares $C_{i,j}=(i-1,i)\times(j-1,j)$, $1\le i,j\le n$.
A line \emph{pierces} a cell if it meets its interior. A family of lines \emph{pierces the board} if every cell is pierced.
Let $p_n$ be the least size of a piercing family.

The problem is Keller's Problem E2862 \cite{Keller}. Ackerman and Pinchasi \cite{AP} proved $p_n\ge\lceil\frac23n\rceil$.
Ambrus, B\'ar\'any, Frankl and Varga \cite{ABFV} proved $0.7n\le p_n\le n-1$ for large $n$ and conjectured $p_n=n-1$.
Alkauskas \cite{Alk26} proved $p_n\ge n-1$ for families of lines whose slopes all lie in $(-\infty,-1]\cup(0,1]$. For unrestricted lines the conjecture is open.
The lower bound of \cite{ABFV} is the method of Theorem \ref{thm:method} below with a quadratic polynomial. We use a polynomial of degree $28$.

For a line $\ell$, let $\gamma$ be the measure on $\ell$ with $d\gamma=|dx|+|dy|$. Explicitly, if $\ell$ is neither horizontal nor vertical,
\[
\int_{\ell\cap[0,1]^2}f\dg=\int f\big(x,\ell(x)\big)\,dx+\int f\big(\ell^{-1}(y),y\big)\,dy ,
\]
the two integrals taken over the projections of $\ell\cap[0,1]^2$ to the $x$-axis and to the $y$-axis. For a horizontal line only the first term is present, and for a vertical line only the second.

\begin{theorem}\label{thm:method}
Let $W$ be a polynomial with $W\ge0$ on $[0,1]^2$. Suppose $\iint_{[0,1]^2}W\,dx\,dy=\alpha$ and $\int_{\ell\cap[0,1]^2}W\dg\le\beta$ for every line $\ell$. Then
\[
\liminf_{n\to\infty}\frac{p_n}n\ \ge\ \frac\alpha\beta .
\]
\end{theorem}

\begin{theorem}\label{thm:main}
There is an explicit polynomial $W_{14}(x,y)$ with $W_{14}\ge0$ on $[0,1]^2$,
\[
\iint_{[0,1]^2}W_{14}\,dx\,dy=\frac{662695081102818845329}{881156433506484375000},\qquad
\int_{\ell\cap[0,1]^2}W_{14}\dg\le\frac{100003}{100000}\ \text{ for every line }\ell .
\]
Consequently
\[
\liminf_{n\to\infty}\frac{p_n}{n}\ \ge\ \frac{2650780324411275381316}{3524731472797958278125}\ =\ 0.75205170\ldots\ >\ \frac34 .
\]
\end{theorem}

Section \ref{sec:method} proves Theorem \ref{thm:method}. Section \ref{sec:main} reduces Theorem \ref{thm:main} to a finite computation, Theorem \ref{thm:computer}, which is carried out by the program described in Section \ref{sec:verifier}.

\section{Proof of Theorem \ref{thm:method}}\label{sec:method}

We scale $B_n$ by $\frac1n$ onto $[0,1]^2$. The scaled cells $C_{i,j}/n$ have side $\frac1n$, and the interior grid lines are $x=\frac in$ and $y=\frac jn$ for $1\le i,j\le n-1$.
For a polynomial $V$ on $[0,1]^2$ let
\[
w_{i,j}=n^2\iint_{C_{i,j}/n}V\,dx\,dy ,
\]
the average of $V$ over the scaled cell, and for a line $\ell$ let
\[
S(\ell)=\sum_{(i,j):\ \ell\ \text{pierces}\ C_{i,j}}w_{i,j} .
\]
We write $\|P-P'\|_1=|x-x'|+|y-y'|$ for the $\ell_1$ distance. If $|\partial_xV|,|\partial_yV|\le\Lambda$ on $[0,1]^2$ then $|V(P)-V(P')|\le\Lambda\,\|P-P'\|_1$ for $P,P'\in[0,1]^2$.

\begin{lemma}\label{lem:discrete}
Let $V$ be a polynomial with $0\le V\le M$ and $|\partial_xV|,|\partial_yV|\le\Lambda$ on $[0,1]^2$. Then for every line $\ell$,
\[
S(\ell)\ \le\ n\int_{\ell\cap[0,1]^2}V\dg\ +\ 3M+8\Lambda .
\]
\end{lemma}

\begin{proof}
Let $\ell^\circ=\ell\cap(0,1)^2$. If $\ell^\circ$ is empty, or $\ell$ is contained in a grid line, then $\ell$ pierces no cell and $S(\ell)=0$. So assume neither.

Call a point of $\ell^\circ$ a \emph{crossing} if it lies on a grid line. Call the crossing \emph{vertical} if its $x$-coordinate is a multiple of $\frac1n$, and \emph{horizontal} otherwise. A horizontal crossing has $y$-coordinate a multiple of $\frac1n$.
There are at most $n-1$ vertical crossings, at most one on each line $x=\frac in$, and at most $n-1$ horizontal ones.
Fix a direction along $\ell$. The crossings cut $\ell^\circ$ into open pieces. Each piece lies in the interior of one cell, and $\ell$ meets each cell interior in at most one segment. So the pieces correspond one to one to the pierced cells.

Every piece except the first begins at a crossing $X$, and $X$ lies in the closure of the piece's cell $C$. So every point of $C$ is within $\ell_1$-distance $\frac2n$ of $X$, and the average of $V$ over $C$ is at most $V(X)+\frac{2\Lambda}n$. The cell of the first piece contributes at most $M$. There are at most $2n$ crossings. Hence
\begin{equation}\label{eq:cross}
S(\ell)\ \le\ M+\sum_{X\ \text{vertical}}V(X)+\sum_{X\ \text{horizontal}}V(X)+4\Lambda .
\end{equation}

If $\ell$ is vertical, then there are no vertical crossings, and \eqref{eq:vert} below holds with the integral equal to $0$.
Otherwise write $y=\ell(x)$, let $t$ be the slope, and let $[x_-,x_+]$ be the projection of $\ell\cap[0,1]^2$ to the $x$-axis. Put $r=x_+-x_-\le1$ and $s=|t|\,r\le1$, the run and the rise of the chord.
The vertical crossings are the points $X_i=(\frac in,\ell(\frac in))$ for the integers $i$ with $x_-<\frac in<x_+$, say $i=a,\dots,b$.
Bound $V(X_a)\le M$. For $a<i\le b$ the part of $\ell$ over $[\frac{i-1}n,\frac in]$ lies in $[0,1]^2$, and each of its points is within $\ell_1$-distance $\frac{1+|t|}n$ of $X_i$. So
\[
V(X_i)\ \le\ n\int_{(i-1)/n}^{i/n}V\big(x,\ell(x)\big)\,dx+\frac{(1+|t|)\Lambda}n .
\]
These $b-a$ intervals are disjoint subintervals of $[x_-,x_+]$, so $b-a\le nr$. Summing, and using $V\ge0$ to extend the integral to $[x_-,x_+]$,
\begin{align}\nonumber
\sum_{X\ \text{vertical}}V(X) &\le\ M+n\int_{x_-}^{x_+}V\big(x,\ell(x)\big)\,dx+r(1+|t|)\Lambda\\
&
\ \le\ M+n\int_{x_-}^{x_+}V\big(x,\ell(x)\big)\,dx+2\Lambda, \label{eq:vert}
\end{align}
since $r(1+|t|)=r+s\le2$.

We now use the same argument with the roles of $x$ and $y$ exchanged. If $\ell$ is horizontal, then there are no horizontal crossings. Otherwise, write $x=\ell^{-1}(y)$, with slope $1/|t|$ in this parametrisation, and let $[y_-,y_+]$ be the projection of $\ell\cap[0,1]^2$ to the $y$-axis. The horizontal crossings are among the points $Y_j=(\ell^{-1}(\frac jn),\frac jn)$ with $y_-<\frac jn<y_+$, and omitting some of these only decreases the sum. The number of nodes after the first is at most $ns$, and the error per node is $\frac{(1+1/|t|)\Lambda}n$. Since $s(1+1/|t|)=s+r\le2$,
\begin{equation}\label{eq:hor}
\sum_{X\ \text{horizontal}}V(X)\ \le\ M+n\int_{y_-}^{y_+}V\big(\ell^{-1}(y),y\big)\,dy+2\Lambda .
\end{equation}

We now use \eqref{eq:vert} and \eqref{eq:hor} to get an upper bound for the right-hand side in the inequality \eqref{eq:cross}. The two integrals add up to $\int_{\ell\cap[0,1]^2}V\dg$ by the definition of $\gamma$ in Section \ref{sec:results}. So $S(\ell)\le n\int_{\ell\cap[0,1]^2}V\dg+3M+8\Lambda$.
\end{proof}

\begin{proof}[Proof of Theorem \ref{thm:method}]
Let $M,\Lambda$ be as in Lemma \ref{lem:discrete} with $V=W$. The scaled cells tile $[0,1]^2$ up to a set of measure zero, so
\[
\sum_{i,j}w_{i,j}=n^2\iint_{[0,1]^2}W\,dx\,dy=\alpha n^2 .
\]
Let $\ell_1,\dots,\ell_{p_n}$ pierce $B_n$. Every cell is pierced by at least one $\ell_r$, and every $w_{i,j}\ge0$. So by Lemma \ref{lem:discrete}
\[
\alpha n^2\ \le\ \sum_{r=1}^{p_n}S(\ell_r)\ \le\ p_n\,(\beta n+3M+8\Lambda).
\]
Divide by $n(\beta n+3M+8\Lambda)$ and let $n\to\infty$.
\end{proof}

\section{Reduction of Theorem \ref{thm:main} to a computation}\label{sec:main}
Recall that the Bernstein polynomials of degree $D$ are
\[
B^D_i(z)=\binom Di z^i(1-z)^{D-i},\qquad i=0,1,\dots,D .
\]
They are non-negative on $[0,1]$, and they form a basis of the polynomials in $z$ of degree at most $D$. The $(D+1)^2$ products $B^D_i(z_1)B^D_j(z_2)$ form a basis of the polynomials in $z_1,z_2$ of degree at most $D$ in each variable.

The polynomial $W_{14}$ of Theorem \ref{thm:main} will be written as
\begin{equation}\label{eq:W14}
W_{14}(x,y)=\sum_{i,j=0}^{14}b_{ij}\,B^{14}_i\big((2x-1)^2\big)\,B^{14}_j\big((2y-1)^2\big).
\end{equation}
It is determined by $225$ rational numbers $b_{ij}$, $0\le i,j\le14$, with denominator $10^6$. These coefficients were found by solving a linear program numerically. The $120$ values with $i\le j$ are listed in the file \texttt{W14.txt} of \cite{Verifier}.

We choose the coefficients so that $b_{ij}=b_{ji}$. Since $W_{14}$ depends on $x$ and $y$ only through $(2x-1)^2$ and $(2y-1)^2$, it is invariant under the symmetries $(x,y)\mapsto(1-x,y)$, $(x,y)\mapsto(x,1-y)$ and $(x,y)\mapsto(1-x,1-y)$ of the square, and since $b_{ij}=b_{ji}$, it is invariant under the symmetry $(x,y)\mapsto(y,x)$.

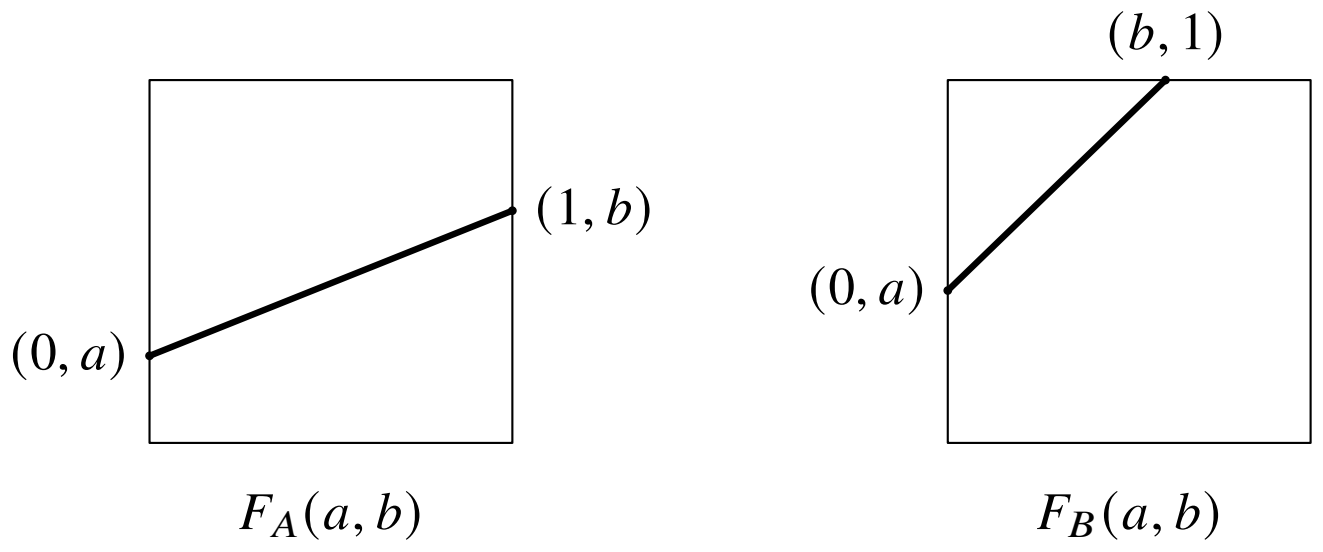
\begin{figure}[ht]
\centering
\begin{tikzpicture}[scale=2.4,line join=round,font=\small]
\begin{scope}
  \draw (0,0) rectangle (1,1);
  \draw[very thick] (0,0.24) -- (1,0.64);
  \fill (0,0.24) circle (0.012) node[left] {$(0,a)$};
  \fill (1,0.64) circle (0.012) node[right] {$(1,b)$};
  \node at (0.5,-0.2) {$F_A(a,b)$};
\end{scope}
\begin{scope}[shift={(2.2,0)}]
  \draw (0,0) rectangle (1,1);
  \draw[very thick] (0,0.42) -- (0.6,1);
  \fill (0,0.42) circle (0.012) node[left] {$(0,a)$};
  \fill (0.6,1) circle (0.012) node[above] {$(b,1)$};
  \node at (0.5,-0.2) {$F_B(a,b)$};
\end{scope}
\end{tikzpicture}
\caption{The segments along which $F_A$ and $F_B$ integrate $W$.}
\label{fig:segments}
\end{figure}

For a polynomial $W$ on $[0,1]^2$ and $(a,b)\in[0,1]^2$ define
\begin{align}
F_A(a,b)&=(1+b-a)\int_0^1W\big(\tau,\ (1-\tau)a+\tau b\big)\,d\tau,\label{eq:FA}\\
F_B(a,b)&=(1-a+b)\int_0^1W\big(b\tau,\ (1-\tau)a+\tau\big)\,d\tau.\label{eq:FB}
\end{align}
Here $F_B(a,b)$ is the integral of $W$ with respect to $\gamma$ along the segment from $(0,a)$ to $(b,1)$, and, when $a\le b$, $F_A(a,b)$ is the integral along the segment from $(0,a)$ to $(1,b)$; see Figure \ref{fig:segments} and the proof below.

\begin{lemma}\label{lem:reduction}
Let $W$ be a polynomial with $W\ge0$ on $[0,1]^2$ that is invariant under the four symmetries $(x,y)\mapsto(1-x,y)$, $(x,y)\mapsto(x,1-y)$, $(x,y)\mapsto(1-x,1-y)$ and $(x,y)\mapsto(y,x)$. Then for every line $\ell$,
\begin{equation}\label{eq:reduction}
\int_{\ell\cap[0,1]^2}W\dg\ \le\ \max\Big(\max_{[0,1]^2}F_A,\ \max_{[0,1]^2}F_B\Big).
\end{equation}
\end{lemma}

\begin{proof}
The four symmetries map lines to lines and preserve $\gamma$ and $W$.
If $\ell\cap[0,1]^2$ is empty or a single point, the integral on the left-hand side of \eqref{eq:reduction} is $0$. The right-hand side is non-negative, because $W\ge0$ and the arguments of $W$ in \eqref{eq:FA} and \eqref{eq:FB} lie in $[0,1]^2$.
Otherwise $\ell\cap[0,1]^2$ is a segment whose endpoints lie on two different edges of the square.

Suppose the endpoints lie on the left and the right edge, say at $(0,a)$ and $(1,b)$. Applying the symmetry $(x,y)\mapsto(x,1-y)$ if necessary, we may assume $a\le b$. The points of the segment are $\big(\tau,(1-\tau)a+\tau b\big)$, $\tau\in[0,1]$, so $dx=d\tau$, $dy=(b-a)\,d\tau$, and $d\gamma=(1+b-a)\,d\tau$. Hence the left-hand side of \eqref{eq:reduction} is $F_A(a,b)$.
If the endpoints lie on the bottom and the top edge, the symmetry $(x,y)\mapsto(y,x)$ reduces this to the previous case.

Suppose the endpoints lie on two adjacent edges. One of the symmetries $(x,y)\mapsto(1-x,y)$, $(x,y)\mapsto(x,1-y)$, $(x,y)\mapsto(1-x,1-y)$, or no symmetry, maps the segment onto a segment from the left edge to the top edge, say from $(0,a)$ to $(b,1)$. Its points are $\big(b\tau,(1-\tau)a+\tau\big)$, $\tau\in[0,1]$, so $dx=b\,d\tau$, $dy=(1-a)\,d\tau$, and $d\gamma=(1-a+b)\,d\tau$. Hence the left-hand side of \eqref{eq:reduction} is $F_B(a,b)$.
\end{proof}

Denote by $F_A^{14}$ and $F_B^{14}$ the polynomials \eqref{eq:FA} and \eqref{eq:FB} for $W=W_{14}$.

\begin{theorem}\label{thm:computer}
The polynomial $W_{14}$ has the following properties.
\begin{enumerate}
\item[(a)] $W_{14}\ge0$ on $[0,1]^2$.
\item[(b)] $\iint_{[0,1]^2}W_{14}\,dx\,dy=\frac{662695081102818845329}{881156433506484375000}$.
\item[(c)] $F_A^{14}\le\frac{100003}{100000}$ and $F_B^{14}\le\frac{100003}{100000}$ on $[0,1]^2$.
\end{enumerate}
\end{theorem}

Theorem \ref{thm:computer} is proved by computer in Section \ref{sec:verifier}.

\begin{proof}[Proof of Theorem \ref{thm:main}]
By Theorem \ref{thm:computer}(a), Lemma \ref{lem:reduction} applies to $W_{14}$, and with Theorem \ref{thm:computer}(c) it gives $\int_{\ell\cap[0,1]^2}W_{14}\dg\le\frac{100003}{100000}$ for every line $\ell$. Together with Theorem \ref{thm:computer}(a) and (b) these are the hypotheses of Theorem \ref{thm:method} with $\beta=\frac{100003}{100000}$, and $\alpha/\beta$ is the fraction in Theorem \ref{thm:main}.
\end{proof}

\section{The verifier}\label{sec:verifier}

Theorem \ref{thm:computer} is proved by the program \cite{Verifier}. It is written in C++, and all its arithmetic is exact. Integers are held in fixed-width words of $1024$ bits from the Boost Multiprecision library, and every operation is checked for overflow. An overflow stops the verification.
The program reads the $120$ coefficients $b_{ij}$ with $i\le j$, sets $b_{ji}=b_{ij}$, forms $W_{14}$ by \eqref{eq:W14}, and verifies the following three statements, where $\beta=\frac{100003}{100000}$:
\begin{align}
&W_{14}\ge0\ \text{ on }[0,1]^2, \tag{V1}\label{V1}\\
&\iint_{[0,1]^2}W_{14}\,dx\,dy=\tfrac{662695081102818845329}{881156433506484375000}, \tag{V2}\label{V2}\\
&\beta-F_A^{14}\ge0\ \text{ and }\ \beta-F_B^{14}\ge0\ \text{ on }[0,1]^2. \tag{V3}\label{V3}
\end{align}
These are parts (a), (b) and (c) of Theorem \ref{thm:computer}.

Put $S=(2x-1)^2$ and $U=(2y-1)^2$. By \eqref{eq:W14}, $W_{14}$ is a polynomial in $(S,U)$ of degree $14$ in each variable, and the map $(x,y)\mapsto(S,U)$ sends $[0,1]^2$ onto itself. So \eqref{V1} says that this polynomial is $\ge0$ on $[0,1]^2$.

The functions $F_A^{14}$ and $F_B^{14}$ are defined by integrals, one for each $(a,b)$. Substituting the linear expressions of \eqref{eq:FA} and \eqref{eq:FB} into $W_{14}$ gives a polynomial in $a$, $b$ and $\tau$, and integrating each power $\tau^k$ over $[0,1]$ gives $\frac1{k+1}$. Hence $F_A^{14}$ and $F_B^{14}$ are polynomials in $(a,b)$, of degree at most $29$ in each variable: $W_{14}$ has degree $28$ in each of $x$ and $y$, the arguments of $W$ in \eqref{eq:FA} and \eqref{eq:FB} are of degree $1$ in each of $a$ and $b$, and the factors in front add $1$. The program computes them exactly, so \eqref{V3} says that two explicit polynomials are $\ge0$ on $[0,1]^2$.

For \eqref{V2}, the program writes $W_{14}$ in the monomials $S^iU^j$ and uses $\iint_{[0,1]^2}S^iU^j\,dx\,dy=\frac1{(2i+1)(2j+1)}$.

Let $g$ be a polynomial of degree at most $D_1$ in $z_1$ and at most $D_2$ in $z_2$. It can be written as
\begin{align}\label{eq:bern}
g(z_1,z_2)=\sum_{i=0}^{D_1}\sum_{j=0}^{D_2}c_{ij}\,B^{D_1}_i(z_1)\,B^{D_2}_j(z_2).
\end{align}
If all $c_{ij}\ge0$, then $g\ge0$ on $[0,1]^2$, because the Bernstein polynomials are non-negative on $[0,1]$. For $W_{14}$ as a polynomial in $(S,U)$, the coefficients in \eqref{eq:bern} are the $b_{ij}$. They are not all non-negative, so \eqref{V1} is not immediate, and the same holds for the two polynomials in \eqref{V3}.

The program therefore subdivides. It divides the square into four squares of half the side, maps each of them onto $[0,1]^2$ by an affine change of variables, and computes the coefficients \eqref{eq:bern} of $g$ on each. On a square where all coefficients are non-negative, $g\ge0$ is proved. A square where some coefficient is negative is divided again. If $g>0$ on $[0,1]^2$, this terminates, since the coefficients of $g$ on a shrinking square converge to values of $g$. The result is a finite list of squares covering $[0,1]^2$, with non-negative coefficients on each, and hence a proof that $g\ge0$.

For \eqref{V1} the final list has $160$ squares, and for the two polynomials in \eqref{V3} it has $175$ and $448$ squares. The whole run takes less than a second on a laptop.

The maximum of $F_A^{14}$ and $F_B^{14}$ on $[0,1]^2$ exceeds $1$ by about $10^{-5}$, which is why $\beta$ is slightly larger than $1$.

\section*{Acknowledgements}
The authors used Anthropic's Claude models (Claude Fable 5 and 5.1 and earlier Claude Opus models, accessed through claude.ai and Claude Code) between July and November 2026,
in developing the ideas of the proofs, in writing the verifier, and in editing the text. The authors take full responsibility for the correctness of the results.

\end{document}